\documentclass[11pt]{article}
\usepackage{amsfonts,amssymb,amsmath,amsthm,amscd}
\usepackage{graphics,graphicx,rotating}
\usepackage{tikz}
\usetikzlibrary{cd}
\usepackage[T1]{fontenc}
\newtheorem{proposition}{Proposition}[section]

\newtheorem{theorem}[proposition]{Theorem}

\DeclareMathOperator{\st}{\;s.t.\,}
\DeclareMathOperator{\et}{\;and\;}

\begin{document}

\newcommand{\EC}{{\sf EC}}
\newcommand{\AS}{{\sf AS}}
\newcommand{\DD}{{\sf D}}

\newcommand{\1}{{\; \raisebox{3pt}{\circle*{2}}}}
\newcommand{\2}{{\; \raisebox{3pt}{\circle*{2} \hspace{1pt} \circle*{2}}}}

\title{Note about holomorphic maps on a compact Riemann surface}
\author{Lucio Guerra}
\date{}
\maketitle

\begin{abstract}
Note to a paper of M. Tanabe concerning the classical theorem of M. De Franchis and F. Severi. %\hfill (2021, September) 
\end{abstract}

\section{Introduction}

\subsection{A classical result}
Let $X$ be a compact Riemann surface, i.e. a projective algebraic curve over the complex field, of genus $g \geq 2$. Let 
$$\mathcal F(X)$$ 
be the collection of holomorphic maps 
$f:X \rightarrow Y$ such that $f$ is surjective and $Y$ is another Riemann surface of genus $\geq 2$.
Two holomorphic maps $f,f'$ from $X$ to $Y,Y'$ are isomorphic if and only if
there is a biholomorphic map $g: Y \rightarrow Y'$ such that $f' = g \circ f$.
$$\begin{tikzcd}
& X \arrow[ld, "f" '] \arrow[rd, "f'"] \\ 
Y \arrow[rr, "g" '] && Y' 
\end{tikzcd}$$
In this case we write $f \simeq f'$. Denote with the symbol
$$\overline{\mathcal F(X)}$$ 
the collection of isomorphism classes $[f]$ of maps $f \in \mathcal F(X)$. 

This is a finite set by a famous theorem of 
M. De Franchis and F. Severi. There is some effective upper bound for the cardinality of this set, i.e. an explicit function $F(g)$, only depending on the genus of the curve, such that $\# \overline{\mathcal F(X)} \leq F(g)$. As far as we know, the best such result is due to M. Tanabe \cite{T}, and in particular the asymptotic behavior of the bounding function is $O((cg)^{5g})$ where $c$ is some `universal' constant, independent of $g$. This work has been extended to higher dimensions by J. C. Naranjo and G. P. Pirola \cite{NP}, for the collection of rational maps of finite degree between two fixed algebraic varieties of general type.

One major idea in the work of Tanabe is that of adding structure to holomorphic maps, in such a way that a kind of classifying correspondence arises.
The aim of this note is to explain in some detail the way how we understand this main idea, using the approach that is suggested in the work of Naranjo and Pirola.

\subsection{Additional structure for holomorphic maps}

Let $f: X \rightarrow Y$ be an element of $\mathcal F(X)$. If

\begin{itemize}
\item[]
$\omega \in H^{1,0}(Y)$ is a holomorphic form, with $\omega \neq 0$, and \smallskip
\item[]
$y \in Y$ is a point in the support of $div(\omega)$, \smallskip
\end{itemize}
we say that $\omega$ as well as
$(\omega, y)$ is a kind of additional structure or marking for $f$,
and that $(f,\omega)$ as well as 
$(f, \omega, y)$ is a marked map. Consider the following collections of maps with additional structure: \medskip \\
$\begin{array}{rcl}
\mathcal F^{\1} (X) & := & \{ (f, \omega) \st f \in \mathcal F(X) \et \omega \in H^{1,0}(Y) \st \omega \neq 0 \}, \medskip \\
\mathcal F^{\2} (X) & := & \{ (f, \omega, y) \st (f, \omega) \in \mathcal F^\1(X) \et y \in | div(\omega) | \}. \end{array}$
\medskip \\
On each collection there is a natural notion of isomorphism, that we denote again with the symbol $\simeq$. Namely, 
$(f, \omega) \simeq (f', \omega')$ if and only if there is an isomorphism $g: Y \rightarrow Y'$ such that $f' = g \circ f$, as in the diagram above, and $\omega = g^\ast (\omega')$; and, moreover, $(f, \omega, y) \simeq (f', \omega', y')$ if and only if in addition $g(y) = y'$.  Denote with the symbols
$$\overline{\mathcal F^{\1}(X)} \et \overline{\mathcal F^{\2}(X)}$$
the collections of isomorphism classes $[(f,\omega)]$ and $[(f,\omega,y)]$, respectively. The natural forgetful maps 
$$\overline{\mathcal F^{\2}(X)} \longrightarrow \overline{\mathcal F^{\1}(X)} \longrightarrow \overline{\mathcal F(X)}$$ are surjective. 

\subsection{Bounding maps with one and the same marking}

Consider the commutative diagram
$$\begin{tikzcd}
\overline{\mathcal F^{\2}(X)} \arrow[d] \arrow[r] & H^{1,0}(X) \times Div^+(X) \arrow[d] \\
\overline{\mathcal F^{\1}(X)} \arrow[r] & H^{1,0}(X) \times \mathbb Z^+
\end{tikzcd}$$
where the arrows act as
$$\begin{tikzcd}
\mbox{$[(f, \omega, y)]$} \arrow[d, mapsto] \arrow[r, mapsto] & (f^\ast(\omega), f^\ast(y)) \\
\mbox{$[(f, \omega)]$} \arrow[r, mapsto] & (f^\ast(\omega), \deg(f))
\end{tikzcd} \hspace{5pt}
\begin{tikzcd}
(\theta, D) \arrow[d, mapsto] \\
(\theta, \deg(D)).
\end{tikzcd}$$

The horizontal arrows in the diagram may be thought of as `classifying correspondences' for marked maps, and 
the first main point in \cite{T} is to provide an upper bound for the cardinality of fibres of the two horizontal arrows. From the quoted paper (cf. Lemma 3 and its proof) we extract the following statement, with a slight improvement.

\begin{theorem} \label{theorem}
Let $(f, \omega, y)$ be an element of $\mathcal F^{\2}(X)$ with 
$d = \deg(f)$. Then
\begin{itemize}
\item[($i$)]
the fibre in $\overline{\mathcal F^{\1}(X)}$ containing $[(f, \omega)]$ has cardinality $$\leq \binom{2g-2}{d} (2g-1)^{d-1};$$
\item[($ii$)]
the fibre in $\overline{\mathcal F^{\2}(X)}$ containing $[(f, \omega, y)]$ has cardinality $$\leq (n+1)^{d-1},$$ where $n = ord_y(\omega)$. \end{itemize}
\end{theorem}

First of all, observe that $(i)$ follows from $(ii)$. This is easy to check by just playing in the diagram above. We only remark that 
the map $[(f, \omega, y)] \mapsto (f^\ast(\omega),f^\ast(y))$ lands into the subset 
$$\{ (\theta, D) \in  H^{1,0}(X) \times Div(X) \st div(\theta) \geq D > 0 \}$$ and that for the restriction of the map $(\theta,D) \mapsto (\theta, \deg(D))$ the fibre over $(\theta,d)$ has cardinality $\binom{2g-2}{d}$. 
So let us proceed to the proof of $(ii)$. 

\section{Preliminary remark on comparing two maps}
\label{preliminary}

Let $f: X \rightarrow Y$ be a holomorphic map, of degree $d$. There is an associated map $Y \rightarrow Div^+_d(X)$ such that $y \mapsto f^\ast(y)$. After the natural identification $Div^+_d(X) \longleftrightarrow X^{(d)}$ with the $d$th symmetric product, the induced map $F: Y \rightarrow X^{(d)}$ is then a holomorphic map. Let $f': X \rightarrow Y'$ be a second holomorphic map, of the same degree $d$, and let $F': Y' \rightarrow X^{(d)}$ be the induced holomorphic map. Later on we will make use of the following remark.

If $f$ and $f'$ have infinitely many common fibres, i.e. if there are infinitely many pairs $(y,y') \in Y \times Y'$ such that $f^\ast(y) = {f'}^\ast(y')$, then there is a biholomorphic isomorphism $g: Y \rightarrow Y'$ such that $f' = g \circ f$.

The proof is easy. Just consider the induced maps
$$\begin{tikzcd}
Y \arrow[r, "F"] & X^{(d)} &
Y' \arrow[l, "F'" '] 
\end{tikzcd}.$$
Because of the hypothesis, one has $F(Y) = F'(Y')$ in $X^{(d)}$ and, since $F$ and $F'$ are injective, this gives the isomorphism $g$.

\section{Comparing two maps with additional structure}

\subsection{Equivalence of marked maps} 

Let $(f, \omega, y)$ and $(f', \omega', y')$ be elements of $\mathcal F^{\2}(X)$. 
We say that they are equivalent and write $(f, \omega, y) \approx (f', \omega', y')$ if and only if
$$f^\ast(\omega) = {f'}^\ast(\omega') \et f^\ast(y) = {f'}^\ast(y').$$ 
Clearly isomorphism $\simeq$ implies equivalence $\approx$.

\subsection{The `ratio' of a pair of maps with equivalent structures} 

Assume that $(f, \omega, y) \approx (f', \omega', y')$. There is $x \in X$ such that $f(x) = y$ and $f'(x) = y'$. From the basic formula
$ord_xf^\ast(\omega) + 1 = (ord_{f(x)}\omega + 1) \deg_x(f)$ and from the hypothesis
$\deg_x(f) =  \deg_x(f')$ it follows that $$ord_y(\omega) = ord_{y'}(\omega') =: n.$$ 

Let $V \et V'$ be coordinate neighborhoods of $y \et y'$ in $Y \et Y'$, so that the points have coordinate $0$, and such that the restricted forms 
$$\omega|_V = {\rm d} k(t) \et \omega'|_{V'} = {\rm d} k'(t')$$
are exact differentials, with 
$$ord_0(k) = ord_0(k') = n+1.$$ 
From the basic local theory of holomorphic functions we have that
$$k(t) = (w(t))^{n+1} \et k'(t') = (w'(t'))^{n+1}$$ 
with $w$ and $w'$ invertible in a neighborhood of $0$, 
each determined up to an $(n+1)$-th root of unity, an element of the group $\mathbb U_{n+1}$. If $w \et w'$ are chosen, shrinking $V \et V'$ if necessary, we may assume that there is a biholomorphic transformation $t' = \gamma(t)$ such that $$w(t) = w'(\gamma(t))$$
and so  there is a biholomorphic map $$g: V \longrightarrow V',$$ which in coordinates is expressed by $\gamma$, such that 
$\omega|_V = g^\ast(\omega'|_{V'})$.

Now write the common fibre as $$f^\ast(y) = {f'}^\ast(y') = x_1 + \ldots + x_d.$$ For each $i=1,\ldots,d$
there is a coordinate neighborhood $U_i$ of $x_i$ in $X$ such that $f(U_i) \subseteq V$ and $f'(U_i) \subseteq V'$ and we may assume that if there is a repetition $x_i = x_j$ then the choice is $U_i = U_j$.
$$\begin{tikzcd}
& U_i \arrow[dl, "f|_{U_i}" '] \arrow[dr, "f'|_{U_i}"] & \\
V && V' 
\end{tikzcd}$$
In terms of local coordinates we write 
$$f|_{U_i}: t = \phi_i(s) \et f'|_{U_i}: t' = \phi'_i(s).$$
The condition $f^\ast(\omega) = {f'}^\ast(\omega')$ means that on $U_i$
$$w(\phi_i(s))^{n+1} = w'(\phi'_i(s))^{n+1}$$
holds and therefore for some $c_i \in \mathbb U_{n+1}$ we must have 
$$c_i\, w(\phi_i(s)) = w'(\phi'_i(s)).$$
Thus we have a sequence $(c_1, \ldots, c_d) \in (\mathbb U_{n+1})^d$. Changing the choice of $w$ or $w'$ implies multiplying $(c_1, \ldots, c_d)$ by some $c \in \mathbb U_{n+1}$. 
Hence we have an element $$[(c_1, \ldots, c_d)] \in (\mathbb U_{n+1})^d / \mathbb U_{n+1}$$ attached to the given pair of equivalent marked maps.

More precisely, here we assume that for each equivalence class of marked maps a choice of the listing $x_1,\ldots,x_d$ (with repetitions according to local degrees) of the common fibre of all marked maps in the equivalence class is assigned. Thus the definition above is well posed.

The element $[(c_1, \ldots, c_d)]$ represents the relation, or the `ratio', existing between two maps which have equivalent additional structures, and we introduce the following amusing notation
$$\frac{(f', \omega', y')}{(f, \omega, y)} := [(c_1, \ldots, c_d)]$$
in view of what is going to follow.

\subsection{Algebraic properties of the ratio}

The following properties are easily checked. \medskip 
\begin{itemize}
\item[$-$]
For two isomorphic $(f, \omega, y) \simeq (f', \omega', y')$ the \medskip ratio is $[(1,\ldots,1)]$.  
\item[$-$]
Given three marked maps $(f, \omega, y)$ and $(f', \omega', y')$ and $(f'', \omega'', y'')$ in the same equivalence class, if the ratio of the first two is $[(c_1, \ldots, c_d)]$ and the ratio of the last two is $[(c'_1, \ldots, c'_d)]$  then the ratio of the first and the last is $[(c_1c'_1, \ldots, c_dc'_d)]$. This can be written as
$$\frac{(f', \omega', y')}{(f, \omega, y)}\; \frac{(f'', \omega'', y'')}{(f', \omega', y')} = \frac{(f'', \omega'', y'')}{(f, \omega, y)}.$$ \medskip
\item[$-$]
The element $[(c_1, \ldots, c_d)]$ only depends on the isomorphism classes $[(f, \omega, y)]$ and $[(f', \omega', y')]$. So we can write
$$\frac{(f', \omega', y')}{(f, \omega, y)} =: \frac{[(f', \omega', y')]}{[(f, \omega, y)]}.$$
\end{itemize} \medskip
Note that the last property follows immediately from the other two.
The next one is the fundamental property of the construction.
\begin{itemize}
\item[$-$]
If $(f, \omega, y) \approx (f', \omega', y')$ are two equivalent marked maps, 
their ratio is $[(1,\ldots,1)]$ if and only if there is isomorphism $(f, \omega, y) \simeq (f', \omega', y')$.
\end{itemize}

The condition on the ratio means that for a suitable choice of coordinates $w$ and $w'$ near $y$ and $y'$ one has $c_i = 1$ for $i=1,\ldots,d$.

If the condition is satisfied, and if we define $g$ and its analytic expression $\gamma$ so that $w = w' \circ \gamma$ as above, then the equality $w(\phi_i(s)) = w'(\phi'_i(s))$ implies that $\phi'_i(s) = \gamma(\phi_i(s))$ holds on $U_i$.
$$\begin{tikzcd}
& U_i \arrow[dl, "f|_{U_i}" '] \arrow[dr, "f'|_{U_i}"] & \\
V \arrow[rr, "g" '] && V'
\end{tikzcd}$$

Define $A := \bigcup_{i=1}^d U_i$. It is an open neighborhood of the common fibre $\{ x_1,\ldots,x_d \}$ in $X$, and we have $g \circ f|_A = f'|_A$. 
Then make smaller both $V$ and $V'$ so that ${f}^{-1}(V)$ and ${f'}^{-1}(V')$ are contained in $A$ and even smaller so that moreover $h(V) = V'$; it follows that $${f}^{-1}(V) = {f'}^{-1}(V') =: U$$ and thus we have the commutative diagram
$$\begin{tikzcd}
& U \arrow[dl, "f|_U" '] \arrow[dr, "f'|_U"] & \\
V \arrow[rr, "g" '] && V'
\end{tikzcd}$$

Finally from the remark in \S \ref{preliminary} we have that the isomorphism $V \rightarrow V'$ extends to an isomorphism $Y \rightarrow Y'$ which gives an isomorphism $(f, \omega, y) \simeq (f', \omega', y')$.

\section{Proof of assertion $(ii)$ in Theorem \ref{theorem}} 

With this preparation, let us come back to the problem stated in the beginning.
Given $(f, \omega, y) \in \mathcal F^{\2}(X)$ with $\deg(f) =: d$, we want to study the set 
$$\{ [(f', \omega', y')] \in \overline{\mathcal F^{\2}(X)} \st f^\ast(\omega) = {f'}^\ast(\omega') \et f^\ast(y) = {f'}^\ast(y') \}.$$ 
This is the object of assertion $(ii)$, which claims that $(n+1)^{d-1}$ is an upper bound for the cardinality of the set.

There is the correspondence from this set to the quotient set $(\mathbb U_{n+1})^d / \mathbb U_{n+1}$
that is given by
$$[(f', \omega', y')] \longmapsto [(c_1,\ldots,c_d)] = \frac{[(f', \omega', y')]}{[(f, \omega, y)]}.$$
The claim is now that this is an injective correspondence, and
this implies assertion $(ii)$.

Let $(f',\omega',y')$ and $(f'',\omega'',y'')$ be elements of $\mathcal F^{\2}(X)$ such that 
$$f^\ast(\omega) = {f'}^\ast(\omega') = {f''}^\ast(\omega'') \et f^\ast(y) = {f'}^\ast(y') = {f''}^\ast(y'')$$  
and moreover
$$\frac{[(f', \omega', y')]}{[(f, \omega, y)]} = \frac{[(f'', \omega'', y'')]}{[(f, \omega, y)]}.$$
From the general formula
$$\frac{[(f', \omega', y')]}{[(f, \omega, y)]}\; \frac{[(f'', \omega'', y'')]}{[(f', \omega', y')]} = \frac{[(f'', \omega'', y'')]}{[(f, \omega, y)]}$$ we obtain
$$\frac{[(f'', \omega'', y'')]}{[(f', \omega', y')]} = [(1,\ldots,1)]$$
and from the last property in the previous section we then have $$[(f', \omega', y')] = [(f'', \omega'', y'')]$$
that completes the proof.

\vfill \noindent
email: {\tt lucioguerra56@gmail.com}

\end{document}